\documentclass[fleqn]{mat01}
\usepackage{times,mathtimy,amssymb,latexsym}
\begin{document}

\setcounter{page}{327}
\firstpage{327}

\def\d{\hbox{d}}

\def\coth{\mbox{\rm coth}}

\newtheorem{theore}{Theorem}
\newtheorem{theor}[theore]{\bf Theorem}
\newtheorem{definit}[theore]{\rm DEFINITION}
\newtheorem{lem}[theore]{Lemma}
\newtheorem{rem}[theore]{Remark}
\newtheorem{propo}[theore]{\rm PROPOSITION}
\newtheorem{coro}[theore]{\rm COROLLARY}

\def\remar{\trivlist \item[\hskip \labelsep{\it Remarks.}]}

\title{The Socle and finite dimensionality of some Banach
algebras}

\markboth{Ali Ghaffari and  Ali Reza Medghalchi}{The Socle and
finite dimensionality of some Banach algebras}

\author{ALI GHAFFARI  and ALI REZA MEDGHALCHI}

\address{Department of Mathematics, Semnan University, Semnan,
Iran\\
\noindent E-mail: ghaffari1380@yahoo.com}

\volume{115}

\mon{August}

\parts{3}

\pubyear{2005}

\Date{MS received 13 September 2004; revised 5 April 2005}

\begin{abstract}
The purpose of this note is to describe some algebraic conditions
on a Banach algebra which force it to be finite dimensional. One
of the main results in Theorem~2 which states that for a locally
compact group $G$, $G$ is compact if there exists a measure $\mu$
in $\hbox{Soc}(L^{1}(G))$ such that $\mu(G) \neq 0$. We also prove
that $G$ is finite if $\hbox{Soc}(M(G))$ is closed and every nonzero
left ideal in $M(G)$ contains a minimal left ideal.
\end{abstract}


\maketitle

\section{Introduction}

Let $A$ be a Banach algebra. The first Arens multiplication on
$A^{**}$ is defined in three steps as follows.

For $a, b$ in $A$, $f$ in $A^{*}$ and $F, G$ in $A^{**}$, the
elements $fa$, $Ff$ of $A^{*}$ and $GF$ of $A^{**}$ are defined by
\begin{equation*}
\langle fa, b \rangle = \langle f, ab \rangle, \quad \langle Ff, a
\rangle = \langle F, fa \rangle, \quad \langle GF, f \rangle =
\langle G, Ff \rangle.
\end{equation*}
We know that $A^{**}$ is a Banach algebra with Arens
multiplication. If $A$ has minimal left ideals, the smallest left
ideal containing all of them is called the left Socle of $A$ and
is denoted by $\hbox{Soc}(A)$. If $A$ does not have minimal left
ideals, we define $\hbox{Soc}(A) = (0)$.

Let $G$ be a locally compact group, $L^{1}(G)$ be its group
algebra, and $M(G)$ be its usual measure algebra. Let
$\hbox{LUC}(G)$ denote the closed subspace of bounded left
uniformly continuous function on $G$, i.e. all $f \in C_{b}(G)$
such that the map $x \rightarrow L_{x}f$ from $G$ into $C_{b}(G)$
is continuous. We know that $L^{1}(G)$ is semisimple and minimal
ideals in $L^{1}(G)$ are generated by minimal idempotent \cite{4}.
Filali \cite{6,7} has studied all the finite dimensional left
ideals of $L^{1}(G), L^{1}(G)^{**}$ and $\hbox{LUC}(G)^{*}$. He
has shown that such ideals exist in $L^{1}(G)$ and $M(G)$ if and
only if $G$ is compact. Baker and Filali \cite{2} proved that
minimal left ideals can be of infinite dimension, and that
compactness of $G$ is not necessary for these ideals to exist in
$L^{1}(G)$ and $M(G)$. For a locally compact abelian group $G$,
Filali \cite{8} has shown that $G$ is compact if and only if
$M(G)$ has minimal ideals. In this paper we will show, among other
things, that if there exists a measure $\mu$ in
$\hbox{Soc}(L^{1}(G))$ such that $\mu(G)\neq 0$, then $G$ is
compact ($G$ is an arbitrary locally compact group). Also some
conditions which are equivalent to finite dimensionality for a
Banach algebra $A$ are given.

\section{Main results}

In this section we will study a Banach algebra $A$ when we set
some conditions on $\hbox{Soc}(A)$ and $\hbox{Soc}(A^{**})$. We
know that $\hbox{Soc}(A)$ has been studied in
\cite{3,4,5,12}.

\begin{propo}$\left.\right.$\vspace{.5pc}

\noindent Let $B$ be an ideal in $A$ such that for any $a \in A${\rm ,}
$Ba = (0)$ implies $a = 0$. Then ${\rm Soc}(A) = {\rm Soc}(B)$.
\end{propo}

\begin{proof} Let $a \in A$ and $Aa$ be a minimal left ideal in
$A$. Since $B$ is a left ideal in $A$, $Ba$ is a left ideal in
$A$. By assumption, $Ba = Aa$. It is easy to see that $Ba$ is a
minimal left ideal of $B$. It follows that $\hbox{Soc}(A)
\subseteq \hbox{Soc}(B)$.

To prove the converse, let $b\in B$ and $Bb$ be a minimal left
ideal in $B$. Since $Bb \neq 0$, there exists $b_{1}\in B$ such
that $b_{1}b\neq 0$. We have $Bb_{1}b\subseteq Ab_{1}b\subseteq
Bb$. By assumption $Bb_{1} b \neq 0$, and so
\begin{equation*}
Bb_{1} b = Ab_{1}b = Bb.
\end{equation*}
If $a \in A$ and $ab_{1}b \neq 0$, then $Bab_{1}b \subseteq
Aab_{1}b\subseteq Ab_{1}b \subseteq Bb$ and $Bab_{1}b\neq 0$. But
$Bb$ is a minimal left ideal in $B$, hence
\begin{equation*}
Aab_{1}b = Ab_{1}b = Bb.
\end{equation*}
Therefore $Ab_{1}b = Bb$ is a minimal left ideal in $A$, which
proves $\hbox{Soc}(B) \subseteq \hbox{Soc}(A)$. Consequently
$\hbox{Soc}(A) = \hbox{Soc}(B)$.
\end{proof}

\setcounter{theore}{0}

\begin{coro}$\left.\right.$\vspace{.5pc}

\noindent Let $G$ be a locally compact group. Then ${\rm Soc}
(L^{1}(G)) = {\rm Soc}(M(G))$.
\end{coro}

\begin{proof}
It is known that $L^{1}(G)$ has a bounded approximate identity
bounded by 1. Let $(e_{\alpha})$ be a bounded approximate identity
in $L^{1}(G)$. Let $\mu \in M(G)$ and $L^{1}(G)\mu = 0$. Since
$C_{0}(G) \subseteq \hbox{LUC}(G) \subseteq L^{\infty}(G)L^{1}(G)$
(see \cite{9,10,11}), for $\psi \in C_{0}(G)$, there exists $f\in
L^{\infty}(G)$ and $\nu \in L^{1}(G)$ such that $\psi = f\nu$. We
have
\begin{align*}
\langle \mu, \psi\rangle &= \langle \mu, f\nu \rangle = \lim
\langle \mu, f\nu * e_{\alpha} \rangle\\
&= \lim \langle e_{\alpha} * \mu, f\nu \rangle = \lim\langle e_{a}
* \mu, \psi \rangle = 0.
\end{align*}
It follows that $\mu = 0$. By Proposition~1, we have
$\hbox{Soc}(L^{1}(G)) = \hbox{Soc}(M(G))$.

In the following Theorem we will provide some conditions on $A$
and $A^{**}$ that are sufficient to guarantee finite
dimensionality.
\end{proof}

\setcounter{theore}{0}

\begin{theor}[\!] Let $A$ be a Banach algebra with a bounded
approximate identity. If ${\rm Soc}(A^{**}) = A^{**}${\rm ,} then
$A$ is finite dimensional.
\end{theor}

\begin{proof} Let $(e_{\alpha})$ be a bounded approximate
identity in $A$ and $E \in A^{**}$ be a weak$^{*}$ limit of a
subnet $(e_{\beta})$ of $(e_{\alpha})$. Then $E$ is a right
identity for $A^{**}$. For a nonzero ideal $\cal{J}$ of $A^{**}$,
we take
\begin{equation*}
\Omega = \{\cal{K}; \cal{K} \ \hbox{is a left ideal in} \ A^{**} \
\hbox{and} \ \cal{J} \ \cap \cal{K} = (0)\}.
\end{equation*}
Let $\cal{M}$ be a maximal element in $\Omega$. Now, let $\cal{I}$
be a minimal left ideal of $A^{**}$. If $\cal{I} \cap (\cal{M} +
\cal{J}) = (0)$, then $(\cal{I} + \cal{M})\cap \cal{J} = (0)$ and
so $\cal{I} + \cal{M} \in \Omega$, i.e, $\cal{M} + \cal{I} =
\cal{M}$. It follows that $\cal{I} \subseteq \cal{M} + \cal{J}$.
If $\cal{I}\cap (\cal{M} + \cal{J}) \neq (0)$, then
$\cal{I}\cap(\cal{M} + \cal{J}) = \cal{I}$ and so $\cal{I}
\subseteq \cal{M} + \cal{J}$. This shows that every minimal left
ideal $\cal{I}$ must be contained in $\cal{M} \oplus \cal{J}$.
Since $\hbox{Soc}(A^{**}) = A^{**}$, we have $\cal{M}\oplus
\cal{J} = A^{**}$. For some $J \in \cal{J}$ and $M \in \cal{M}$,
we can write $E = M + J$. It follows that $\cal{J}^{2} = \cal{J}$,
i.e. $\cal{J} \neq (0)$. This shows that $A^{**}$ is semiprime. By
Theorem~5 of \cite{13}, $A^{**}$ is finite dimensional and so $A$
is finite dimensional.
\end{proof}

\begin{coro} $\left.\right.$\vspace{.5pc}

\noindent Let $G$ be a locally compact group. Then $G$ is finite if and only
if~~${\rm Soc}(L^{1}(G)^{**}) = L^{1}(G)^{**}$.
\end{coro}

\begin{proof}
Since $L^{1}(G)$ has a bounded approximate identity, by Theorem~1,
$L^{1}(G)^{**}$ is finite dimensional. Therefore $G$ is finite.
The converse is clear.

Let $A$ be a Banach algebra and let $\hbox{Comp}(A)$ be the
compactum of $A$, that is the set of all $x$ in $A$ such that the
mapping $a\rightarrow x ax$ is a compact operator of $A$ into
itself. Al-Moajil \cite{1} gives some characterizations of a
finite dimensionality of a semisimple Banach algebra in terms of
its compactum and Socle.

For a locally compact abelian group $G$, Filali \cite{8} has shown
that $G$ is compact if and only if $M(G)$ has minimal ideals. In
the following Theorem we set a condition on $\hbox{Soc}(L^{1}(G))$
and prove that $G$ is compact.
\end{proof}

\setcounter{theore}{1}
\begin{theor}[\!]
Let $G$ be a locally compact group. Then $G$ is compact if any of
the following conditions hold{\rm :}
\begin{enumerate}
\renewcommand\labelenumi{\rm (\arabic{enumi})}
\leftskip .1pc
\item There exists a measure $\mu$ in ${\rm Soc}(L^{1}(G))$ such
that $\mu (G) \neq 0${\rm ;}

\item There exists a measure $\mu$ in ${\rm Soc}(M(G))$ such
that $\mu (G) \neq 0${\rm ;}

\item $G$ is an abelian group and ${\rm Soc}(L^{1}(G))\neq 0$.
\end{enumerate}
\end{theor}

\begin{proof}
Let (1) hold. We assume to the contray that $G$ is not compact.
Let $\Omega$ denote the set of all compact subset of $G$ and we
make $\Omega$ directed by $K_{1} \leq K_{2}$ if and only if
$K_{1}\subseteq K_{2}$ for every $K_{1}$ and $K_{2}$ in $\Omega$. For
every $K \in \Omega$, we can choose $g_{K} \notin K$. Without loss
of generality, we may assume that $\delta_{g_{K}} \rightarrow m$ $(m
\in \hbox{LUC}(G)^{*})$ in the $\sigma (\hbox{LUC}(G)^{*}, \
\hbox{LUC}(G))$-topology. It is easy to see that $\langle m, \psi
\rangle = 0$ for every $\psi \in C_{0}(G)$. Also, we have $\langle
\mu m \mu, \psi \rangle = 0$, for every $\psi \in C_{0}(G)$ and
$\mu \in L^{1}(G)$. (Indeed, the formulas which define the first
Arens product in $L^{1}(G)^{**}$ can be used to define a Banach
algebra structure on $\hbox{LUC}(G)^{*}$.)

Choose $\mu \in\hbox{Soc}(L^{1}(G))$ with $\mu(G)\neq 0$. For a
bounded approximate identity $(e_{\alpha})$ in $L^{1}(G)$ with
norm 1 and $g\in G$, we have
\begin{equation*}
\mu * e_{\alpha} * \delta_{g} * \mu \in \{\mu * \nu * \mu ; \nu
\in L^{1}(G), \| \nu \| \leq 1\}.
\end{equation*}
Therefore
\begin{equation*}
\mu * \delta_{g} * \mu \in cl \{\mu * \nu * \mu ; \nu \in
L^{1}(G), \| \nu \| \leq 1\},
\end{equation*}
where closure is taken in the norm topology. But by Proposition~3
of \cite{1}, the set \hbox{$\{\mu*\nu*\mu$}; $\nu \in L^{1}(G), \| \nu
\| \leq 1\}$ is relatively compact. Hence without loss of
generality we may assume that \hbox{$\mu*\delta_{g_{K}}*\mu
\rightarrow \eta\ (\eta \in L^{1}(G))$} in the norm topology. On the
other hand, $\mu\,*$ \hbox{$\delta_{g_{K}}*\mu \rightarrow 0$} in the
$\sigma (M(G), C_{0}(G))$-topology. It follows that $\eta = 0$.
But \hbox{$\mu*\delta_{g_{K}}\,*$} $\mu (G) = \mu(G)^{2}\neq 0$. This
contradicts the fact that $\eta = 0$. Hence $G$ is compact.

Now, let (2) hold. By Corollary~1, $\hbox{Soc}(L^{1}(G)) =
\hbox{Soc}(M(G))$. By (1), $G$ is compact. Let (3) hold. By
\cite{8}, $G$ is compact.

For a locally compact group $G$, $\hbox{Comp}(L^{1}(G))\subseteq
\hbox{Comp} (M(G))$. Indeed, since $L^{1}(G)$ has a bounded
approximate identity with norm 1, for any $\mu \in L^{1}(G)$ we
have
\begin{equation*}
\{\mu * \nu * \mu; \nu \in M(G), \| \nu\| \leq 1\} \subseteq
cl\{\mu
* \eta
* \mu; \eta \in L^{1}(G), \| \eta \| \leq 1\}.
\end{equation*}
If $G$ is a compact group, then for any $\mu \in L^{1}(G)$ both
mapping $\rho_{\mu}$ and $\lambda_{\mu}$ from $L^{1}(G)$ into
$L^{1}(G)$ are compact, where $\rho_{\mu}(\nu) = \nu * \mu$ and
$\lambda_{\mu}(\nu) = \mu * \nu$ for $\nu \in L^{1}(G)$. It
follows that $L^{1}(G) = \hbox{Comp}(L^{1}(G))\subseteq \
\hbox{Comp}(M(G))$ and so $\hbox{Soc}(M(G)) \neq 0$ (Proposition~3
of \cite{1}). By Corollary~1, $\hbox{Soc}(L^{1}(G))\neq 0$.
\end{proof}

\begin{theor}[\!]
Let $A$ be a semiprime Banach algebra with an identity. If ${\rm
Soc}(A)$ is closed and every nonzero left ideal $I$ of $A$
contains a minimal left ideal{\rm ,} then $A$ is finite dimensional.
\end{theor}

\begin{proof}
If $\hbox{Soc}(A) = A$, then $A$ is finite dimensional (Theorem~5
of \cite{13}). Otherwise we can find a sequence of pairwise
orthogonal idempotents $(e_{n})$ such that $e_{n} \in \
\hbox{Soc}(A)$, since every nonzero left ideal $I$ of $A$ contains
a minimal left ideal. By assumption $\hbox{Soc}(A)$ is closed, so
\begin{equation*}
a = \sum\limits_{n = 1}^{\infty} \frac{e_{n}}{2^{n}\| e_{n} \|}
\in \hbox{Soc} (A).
\end{equation*}
Therefore, the sequence $e_{n}$ is contained in $aAa$, and since
it is an infinite set and linearly independent by the
orthogonality of its elements, we have $aAa$ is infinite
dimensional. This contradicts the fact that $a\in \hbox{Soc}(A)$
(Lemma~2 of \cite{1}). Consequently $\hbox{Soc}(A) = A$, and so
$A$ is finite dimensional.
\end{proof}

\setcounter{theore}{2}
\begin{coro}$\left.\right.$\vspace{.5pc}

\noindent Let $G$ be a locally compact group such that every nonzero ideal
$I$ in $M(G)$ contains a minimal left ideal and let ${\rm Soc}
(M(G))$ be closed. Then $G$ is finite.
\end{coro}

\begin{proof}
See Theorem~3.
\end{proof}


\begin{thebibliography}{99}
\bibitem{1} Al-Moajil~A~H, The compactum and finite dimensionality
in Banach algebras, {\it Int. J. Math. \& Math. Sci.} {\bf 5}
(1982) 275--280

\bibitem{2} Baker~J~W and Filali~M, On minimal ideals in some Banach
algebras associated with a locally compact group, {\it J. London
Math. Soc.} {\bf 63} (2001) 83--98

\bibitem{3} Bresar~M and Semrl~P, Finite rank elements in semisimple
Banach algebras, {\it Studia Math.} {\bf 128} (1998) 287--298

\bibitem{4} Dales~H~G, Banach algebras and automatic continuity (New York, Oxford:
Oxford University Press Inc.) (2000)

\bibitem{5} Dalla~L, Giotopoulos~S and Katseli~N, The Socle and
finite dimensionality of a semiprime Banach algebra, {\it Studia
Math.} {\bf 92} (1989) 201--204

\bibitem{6} Filali~M, Finite dimensional left ideals in some Banach
algebras associated with a locally compact group, {\it Proc. Am.
Math. Soc.} {\bf 127} (1999) 2325--2333

\bibitem{7} Filali~M, Finite dimension right ideals in some Banach
algebras associated with a locally compact group, {\it Proc. Am.
Math. Soc.} {\bf 127} (1999) 1729--1734

\bibitem{8} Filali~M, The ideal structure of some Banach algebras,
{\it Math. Proc. Camb. Philos. Soc.} {\bf 111} (1992) 567--576

\bibitem{9} Ghaffari~A, Convolution operators on semigroup
algebras, {\it Southeast Asian Bull. Math.} {\bf 27} (2004)
1025--1036

\bibitem{10} Hewitt~E and Ross~K~A, Abstract harmonic analysis
(Heidelberg and New York: Springer-Verlag, Berlin) (1963) vol.~1

\bibitem{11} Hewitt~E and Ross~K~A, Abstract harmonic analysis
(Heidelberg and New York: Springer-Verlag, Berlin) (1970) vol.~II

\bibitem{12} Takahasi~S~E, Finite dimensionality in socle of Banach
algebras, {\it Int. J. Math \& Math. Sci.} {\bf 7} (1984) 519--522

\bibitem{13} Tullo~A~W, Condition on Banach algebras which imply
finite dimensionality, {\it Proc. Edinburg Math. Soc.} {\bf 20}
(1976) 1--5
\end{thebibliography}
\end{document}